\documentclass{amsart}
\usepackage{latexsym,amsxtra,amscd,ifthen}
\usepackage{amsfonts}
\usepackage{color,graphicx}
\usepackage{verbatim}
\usepackage{amsmath}
\usepackage{amsthm}
\usepackage{amssymb}
\usepackage{cleveref}
\usepackage{url}
\usepackage{todonotes}
\usepackage{soul}

\usepackage{tikz}
\usetikzlibrary{arrows,shapes,chains}

\theoremstyle{plain}

\newtheorem{theorem}{Theorem}
\newtheorem*{theorem*}{Theorem}

\newtheorem{corollary}[theorem]{Corollary}
\newtheorem{conjecture}[theorem]{Conjecture}

\numberwithin{theorem}{section}
\numberwithin{equation}{theorem}

\theoremstyle{definition}
\newtheorem{definition}[theorem]{Definition}
\newtheorem{example}[theorem]{Example}

\newtheorem{question}[theorem]{Question}
\newtheorem*{question*}{Question}
\usepackage[normalem]{ulem}

\newcommand{\Z}{\mathbb{Z}}
\newcommand{\fm}{\mathfrak{m}}

\DeclareMathOperator{\LND}{LND}
\DeclareMathOperator{\ML}{ML}

\DeclareMathOperator{\kk}{\Bbbk}
\DeclareMathOperator{\Der}{\textrm{Der}}

\begin{document}

\title[ZCP for noncommutative and Poisson algebras]{A survey on Zariski cancellation problems\\for noncommutative and Poisson algebras}
\date{November 2022}

\author[Huang]{Hongdi Huang}
\address{(Huang) Department of Mathematics, Rice University, Houston, TX 77005, U.S.A.}
\email{h237huan@rice.edu}

\author[Tang]{Xin Tang}
\address{(Tang) Department of Mathematics \& Computer Science, 
Fayetteville State University, Fayetteville, NC 28301,
U.S.A.}
\email{xtang@uncfsu.edu}

\author[Wang]{Xingting Wang}
\address{(Wang) Department of Mathematics, Louisiana State University, Baton Rouge, Louisiana 70803, USA}
\email{xingtingwang@math.lsu.edu}

\begin{abstract}
In this article, we discuss some recent developments of the Zariski Cancellation Problem in the setting of noncommutative algebras and Poisson algebras. 
\end{abstract}



%\thanks{ }

\date{\today}
\maketitle

%\tableofcontents

% \setcounter{section}{-1}
\section{Introduction}
\label{xxsec0}
As remarked by Kraft in his 1995 survey \cite{Kr} that ``...\emph{there is no doubt that complex affine $n$-space $\mathbb{A}^n=\mathbb{A}^n_{\mathbb{C}}$ is one of the basic objects in algebraic geometry. It is therefore surprising how little is known about its geometry and its symmetries. Although there has been some remarkable progress in the last few years, many basic problems remain open.}''  

Among these open questions, features the ``Zariski Cancellation Problem". A birational version of the cancellation problem was first raised by Zariski during the 1949 Paris Colloquium on Algebra and Theory of Numbers \cite{Se}. Throughout this paper, we work over a base field $\kk$. 
All the algebras are over $\kk$ and we reserve $\otimes$ to mean $\otimes_{\kk}$ unless stated otherwise. Let $K$ and $K^{\prime}$ be two finitely generated fields over $\kk$. Denote by $K(x)$ and $K^{\prime}(x)$ the simple transcendental extensions of $K$ and $K^{\prime}$, respectively. Zariski asked whether the isomorphism $K(x)\cong K^{\prime}(x)$ implies $K\cong K^{\prime}$? Over the years, a `biregular' version of the cancellation problem has attracted a great deal of research attention. More precisely, let $\mathbb{A}^{n}=\kk^{n}$ be the affine space. One has the following ``Zariski Cancellation Problem" (abbreviated as {\bf ZCP}). 

\begin{question}[Zariski Cancellation Problem]
Does an isomorphism \[Y\times \mathbb{A}^1~\cong~\mathbb{A}^{n+1}\] imply an isomorphism $Y~\cong ~\mathbb{A}^n$, for any affine variety $Y$?
\end{question}
Algebraically, the Zariski Cancellation Problem asks whether an algebra isomorphism $A[x_{1}]\cong \kk[x_{1}, \cdots, x_{n+1}]$ implies that $A$ is isomorphic to $\kk[x_{1}, \cdots, x_{n}]$ as $\kk$-algebras. The Zariski Cancellation Problem has been solved for $A=\kk[x]$ in \cite{AEH} and for $A=\kk[x,y]$ in \cite{Fu, MS, Ru}. Recently, the Zariski Cancellation Problem has been settled negatively in \cite{Gu1, Gu2} for $\kk[x_{1}, \cdots, x_{n}]$ where $n\geq 3$ and ${\rm char}(\kk)>0$. The Zariski Cancellation Problem remains open for $A=\kk[x_{1},\cdots, x_{n}]$ for $n\geq 3$ in the zero characteristic case. The Zariski Cancellation Problem and its analogs have been further studied for commutative domains, see \cite{Cr, CM, Da,  Fi, Fu, Ho, Ma1} for a sample of references. We refer the reader to the excellent survey papers \cite{Gu3, Gu4} for a complete account in the setting of commutative domains.

In this survey article, we would like to discuss some recent developments in the Zariski Cancellation Problem and its generalizations in the setting of noncommutative algebras and Poisson algebras. Note that the study of cancellation properties for noncommutative rings or algebras dates back to the early 1970s, see for example \cite{AEH, As, BR, CE, EH}. There is a vast literature on the cancellation problem for noncommutative algebras and Poisson algebras; see for instance the incomplete list \cite{BHHV, BZ1, BZ2, Be, GaW, GaWY, LeWZ, LuWZ, TVZ, TZZ}.

 In general, it is natural to identify noncommutative algebras $A$ with similar cancellation properties. If $A$ is a finitely generated $\kk$-algebra, then we will simply call $A$ an \emph{affine algebra} over $\kk$. We denote by $Z(A)$ the center of $A$. The following notions are formulated in \cite{BZ1}.
\begin{definition} \cite[Definition 1.1]{BZ1} \label{def:zcp}
Let $A$ be a $\kk$-algebra.
\begin{enumerate}
\item[(1)] $A$ is said to be
{\it cancellative} if any
$\Bbbk$-algebra isomorphism $A[t]\cong B[t]$ for some $\Bbbk$-algebra $B$ implies that $A\cong B$.
\item[(2)] $A$ is said to be {\it strongly cancellative} if, for each $n \geq 1$,  any
$\Bbbk$-algebra isomorphism $A[t_1, \cdots, t_n]\cong B[t_1, \cdots, t_n]$
for some algebra $B$ implies that $A\cong B$.
\item[(3)] $A$ is said to be {\it universally cancellative} if, for every affine commutative $\kk$-domain $R$ such that the natural map $\kk \to R \to R/I$ is an isomorphism for some ideal $I \subset R$ and every $\kk$-algebra $B$, $A \otimes R \cong B \otimes R$ implies that $A \cong B$.
\end{enumerate}
\end{definition}
It follows from the definition that ``$A$ is universally cancellative" implies ``$A$ is strongly cancellative" and ``$A$ is strongly cancellative" implies ``$A$ is cancellative". Moreover, the size of $Z(A)$ seems to play a significant role in the cancellation problem as shown below.
\begin{theorem}\label{thm:center}
\cite[Proposition 1.3]{BZ1} Let $A$ be a $\kk$-algebra with center being $\kk$. Then $A$ is
universally cancellative.
\end{theorem}
This is a useful result since there are many noncommutative algebras that have their centers equal to the base field $\kk$ and are thus universally cancellative. We have the following typical example from \cite{BZ1}. 
We refer the reader to \cite{LV} for the calculation of centers for many classes of noncommutative algebras.
\begin{example}
\cite[Example 1.4]{BZ1}
(1) Let $\kk$ be a field of characteristic zero and $A_{n}$ be the $n$th Weyl algebra over $\kk$. Then $Z(A_{n}) =\kk$. So $A_{n}$ is universally cancellative.\\
(2) Let $\kk$ be a field and $q \in \kk^{\times}$. Let $A=\kk_{q} [x_{1},\cdots, x_{n}]$ be the skew polynomial ring
generated by $x_{1},\cdots, x_{n}$ subject to the relations $x_{j}x_{i}=q x_{i}x_{j}$ for all $1\leq i<j\leq n$. If $n\geq 2$ and $q$ is not a root of unity, then $Z(A) = \kk$. So $A$ is universally
cancellative.
\end{example}
On the other hand, the centers of many PI noncommutative algebras are strictly larger than the base field $\kk$. Thus novel ideas are needed to address the cancellation problem for noncommutative algebras. We will discuss more about the cancellation problem for noncommutative algebras, with a focus on many important families of Artin--Schelter regular algebras in Section $2$, where we will also address the cancellation problem for non-domains. In Section $3$, we will discuss certain generalized cancellation problems such as Morita and skew cancellations. Section $4$ is devoted to the Zariski Cancellation Problem for Poisson algebras. We will employ the notion of the {\it Gelfand--Kirillov (GK) dimension} instead of the classical Krull dimension. We refer the reader to \cite{KL} for the definition and basic properties of the GK-dimension. It is well known that the GK-dimension of any affine commutative $\kk$-algebra is always equal to its Krull dimension.

\subsection*{Acknowledgements} 
The authors thank the referee for carefully reading the paper and providing many thoughtful comments and also thank Milen Yakimov for valuable correspondence on the subject. Part of this paper was formulated during the conference ``Recent Advances and New Directions in the Interplay of Noncommutative Algebra and Geometry" in June 2022 and the BIRS Workshop on ``Noncommutative Geometry and Noncommutative Invariant Theory" in September 2022. The authors want to thank the University of Washington and Banff International Research Station for their hospitality and support. Wang was partially supported by Simons collaboration grant \#688403 and AFOSR grant FA9550-22-1-0272.

\section{Zariski Cancellation for Noncommutative Algebras}
In this section, we will discuss the cancellative properties for noncommutative algebras. Among them is an important family called Artin--Schelter (AS) regular algebras,  which are considered noncommutative analogs of polynomial algebras. We refer the reader to \cite{ArS} for the definition and basic properties of AS-regular algebras. AS-regular algebras should serve as important testing examples for the Zariski Cancellation Problem. 

There are two crucial ingredients in the treatment of ZCP for noncommutative algebras, namely the {\it Makar-Limanov invariant} and the {\it discriminant}. The Makar-Limanov invariant was originally introduced by Makar-Limanov \cite{Ma1} and has become a very useful invariant in commutative algebra. It was first used in \cite{BZ1} to solve rigidity, automorphism, isomorphism, and embedding problems for several classes of noncommutative algebras. In the following three subsections, we state the cancellativity results according to the GK-dimension of an algebra.

\subsection{Definitions}
In this subsection, we list some basic definitions and terminologies. We recall the definition of the Makar-Limanov invariant and its variations from \cite{BZ1}, see also \cite{Ma2}. Let $A$ be a $\kk$-algebra and $\Der(A)$ be the collection of all $\kk$-linear derivations of $A$. Denote by $\LND(A)=\{ \delta \in \Der(A) \mid \delta \text{ is locally nilpotent} \}$. When ${\rm char}(\kk)>0$, one needs to use higher derivations. We denote by ${\rm LND}^{H}(A)$ the set of all locally nilpotent higher derivations of $A$ and ${\rm LND}^{I}(A)$ the set of locally nilpotent iterative higher derivations of $A$. We refer the reader to \cite{BZ1, HS, LL} for more details on higher derivations. We have the following definition of rigidity in terms of derivations. Let $^{\ast}$ be either blank, $^H$, or $^I$.

\begin{definition}
\cite[Definition 2.3]{BZ1} 
Let $A$ be an algebra over $\kk$.
\begin{enumerate}
\item [(1)] The \emph{Makar-Limanov$^{\ast}$ invariant} of $A$ is defined to be 
\[
\ML^{*}(A)=\bigcap_{\delta \in \LND^{\ast}(A)} \ker(\delta)
\]
This means that we have the original ${\rm ML}(A)$, as well as, ${\rm ML}^{H}(A)$ and ${\rm ML}^{I}(A)$.
\item [(2)] We say that $A$ is ${\rm LND}^{\ast}$-rigid if ${\rm ML}^{\ast}(A)=A$, or ${\rm LND}^{\ast}(A)=\{0\}$.
\item [(3)] $A$ is called {\it strongly} ${\rm LND}^{\ast}$-rigid if ${\rm ML}^{\ast}(A[t_{1}, \cdots, t_{n}])=A$, for all $n\geq 1$.
\end{enumerate}
\end{definition}

The Makar-Limanov invariant of an algebra $A$ plays a critical role in controlling the cancellation property for $A$. Indeed, we have the following result from \cite{BZ1}. 

\begin{theorem}\label{thm:LND}
\cite[Theorem 3.3]{BZ1} Let $A$ be an affine $\kk$-domain of finite GK-dimension. When $\ast$ is blank we further assume that $A$ contains $\Z$.
\begin{enumerate}
    \item [(1)] If $A$ is strongly ${\rm LND}^{\ast}$-rigid, then $A$ is strongly cancellative.
    \item [(2)] If ${\rm ML}^{\ast}(A[t])=A$, then $A$ is cancellative. 
\end{enumerate}
\end{theorem}

Therefore, to establish the cancellation property for an algebra $A$, it suffices to prove that the algebra $A$ is ${\rm LND}^{\ast}$-rigid. In order to achieve that goal, another invariant called \emph{discriminant} is introduced for any noncommutative algebra that is module-finite over its center. 

\begin{definition}\cite[Definition 1.3]{CPWZ1}
  Let $A$ be a $\kk$-algebra that is a finitely generated free module over its center $Z$ of rank $w$. Let $\{z_1,\ldots,z_w\}$ be a $Z$-basis of $A$ and let tr \(: A \to Z\) be a trace map. Then the \emph{$w$-discriminant} of $A$ over $Z$ is defined to be 
 \[d_w(A/Z) := \!_{Z^\times} \textnormal{det}(\textnormal{tr}(z_iz_j))_{w\times w}.\]
\end{definition}
It is shown in \cite{BZ1} that certain properties such as being dominating and effective of the discriminant control the ${\rm LND}^{\ast}$-rigidity.  We refer the reader to \cite{BZ1, CPWZ1, CPWZ2} for the definition of the notions of effectiveness and dominance of an element in an algebra and their relations to automorphism groups. 

\begin{example}
 \cite[Example 5.6]{CPWZ1} 
Suppose ${\rm char}(\kk)\neq 2$. Consider the following algebra  
\[A=\kk\langle x,y \rangle/(x^2y - yx^2, xy^2 - y^2x, x^6 - y^2).\]
A computation shows that the center $Z$ of $A$ is the polynomial ring generated by $x^2$ and $xy+yx$ and its $4$-discriminant is given by $d_4(A/Z)=(xy-yx)^4$ which is both dominating and effective.
\end{example}

\subsection{Algebras of GK-dimension one or two}
It is well known that the coordinate ring of an affine curve is cancellative \cite{AEH}. In light of this result, one is led to consider the cancellativity for every affine prime $\kk$-algebra of GK-dimension one. Many nice properties of algebras of GK-dimension one are investigated in \cite{SW}, which are useful to characterize the cancellativity. We remark that cancellativity is prevalent in GK-dimensions 1 and 2, while many specific families of algebras are known to be cancellative in higher dimensions \cite{BZ1, LY}.

For an algebra $A$ of GK-dimension one, some affirmative answers to the Zariski cancellation problem from \cite{LeWZ} and \cite{BHHV} are listed in the following theorem.

\begin{theorem}\label{th:1}
Let $A$ be an affine $\kk$-algebra of GK-dimension one. 
\begin{enumerate}
\item \cite[Theorem 0.6]{LeWZ} 
If $\kk$ is algebraically closed and $A$ is prime, then $A$ is cancellative.

 \item \cite[Theorem 1.1]{BHHV}
 If $\kk$ is of characteristic zero and $A$ is a domain, then $A$ is cancellative.
 \end{enumerate}
\end{theorem}

Notice that in the above theorem, for part (a) the authors invoke Tsen's theorem which requires the base field $\kk$ to be algebraically closed. Part (b) is somewhat orthogonal to part (a) since domains of GK-dimension one over algebraically closed fields are commutative by an application of Tsen's theorem to a result of Small and Warfield \cite{SW}. Hence the only case of part (b) covered by (a) is when the algebra is commutative, which was previously known from the result of Abhyankar-Eakin-Heinzer \cite{AEH}.

We remark that the base field $\kk$ plays an important role here, indeed. As mentioned before, $\kk[x_1, \cdots, x_n]$ with $n\geq 3$ is not cancellative when ${\rm char}(\kk)>0$ by \cite{Gu1,Gu2}. When it comes to a noncommutative $\kk$-algebra $A$, where $\kk$ has positive characteristic and is not algebraically closed, a family of non-cancellative algebras is constructed in \cite{BHHV} as below.
\begin{example}\label{counterex}
Let $p$ be a prime, and let $K=\mathbb{F}_p(x_1,\ldots ,x_{p^2-1})$. Let $k=\mathbb{F}_p(x_1^p,\ldots, x_{p^2-1}^p)$ and we let $\delta$ be the $k$-linear derivation of $K$ given by
$\delta(x_i)=x_{i+1}$ for $i=1,\ldots ,p^2-1$, where we take $x_{p^2}=x_1$.  Since $k$ has characteristic $p>0$,  $\delta^{p^i}$ is a $k$-linear derivation for every $i\ge 0$, and since $\delta^{p^2}(x_i)=\delta(x_i)=x_{i+1}$ for $i=1,\ldots ,p^2-1$, $\delta^{p^{j+2}}=\delta^{p^j}$ for every $j\ge 0$.  We let $\delta':=\delta^p$, which as we have just remarked is a $k$-linear derivation of $K$. We let $A=K[x;\delta]$ and $B=K[x';\delta']$.  Since ${\rm ad}_u^p = {\rm ad}_{u^p}$ for any element $u$ in a ring of characteristic $p$, we have $z:=x^{p^2}-x$ and $z':=(x')^{p^2}-x'$ are central by the above remarks. One can check that $A$ and $B$ have Gelfand-Kirillov dimension one and $A[t]\cong B[t']$, but $A\not\cong B$.
\end{example}

It is worth pointing out that the above examples indicate the necessity of the condition that $\kk$ is algebraically closed in the positive characteristic case. Next, we summarize the results of cancellativity of an algebra having GK-dimension two.

\begin{theorem}
\cite[Theorem 0.5]{BZ1}. Let $\kk$ be an algebraically closed field of characteristic zero. Let
$A$ be an affine $\kk$-domain of GK-dimension two. If $A$ is not
commutative, then $A$ is cancellative.
\end{theorem}
Note that, indeed,  there are examples of commutative affine domains of GK-dimension two that are not cancellative \cite{Da, Fi, Fu, Ho}. We quote some counter-examples from \cite{Da} below. 
\begin{example}\label{ex}
 \begin{enumerate}
\item [(a)]Let \(n \geq 1\) and let \(B_{n}\) be the coordinate ring of the surface
\(x^{n} y=z^{2}-1\) over \(\mathbb{C} .\) Then \(B_{i} \not\cong B_{j}\) if \(i \neq j,\) but \(B_{i}[t] \cong B_{j}[t]\) for
all \(i, j \geq 1 .\) Therefore, all the \(B_{n}\)'s are not cancellative.
\item [(b)] Let $A$ and $B$ be the affine coordinate rings of the surfaces $xy-z^{2}+1=0$ and $x^{2}y-z^{2}+1=0$. We have $A\not\cong B$ while $A[t]\cong B[t]$. 
\end{enumerate} 
 \end{example}

\subsection{Algebras with `small' center or higher GK-dimension}
As shown in Theorem \ref{thm:center}, the center $Z(A)$ of an algebra $A$ provides some control on the cancellation property of $A$. Intuitively, if $A$ has a smaller center, then $A$ is closer to be cancellative. So it is natural to ask if 
the cancellativity of the center can imply that of the algebra itself. 

\begin{question}\cite[Question 0.1]{LeWZ}\cite[Conjecture 2.1]{BHHV}
Let $A$ be an affine noetherian domain. Suppose that $Z(A)$ is affine and cancellative. Is $A$ cancellative? In particular, if the center of $A$ is finite-dimensional over $\kk$, is $A$ cancellative?
\end{question}

In \cite{LeWZ}, cancellation property was proved for non-domain noncommutative algebras using a generalization of ideas from \cite{BZ1} relating Makar-Limanov invariant and {\rm LND}-rigidity to their `small' centers, that is, the ${\rm ML}_Z$-invariant and the notion of ${\rm LND}_Z$-rigidity. 

\begin{theorem} The following algebras $A$ are cancellative:
\begin{enumerate}
\item {\rm(\cite[Theorem 0.2]{LeWZ})}  $A$ is strongly Hopfian with artinian center.
\item {\rm(\cite[Corollary 0.3]{LeWZ})} $A$ is left (or right) artinian; in particular, $A$ is finite-dimensional over a base field $\kk$.
\item {\rm(\cite[Corollary 0.3]{LeWZ})} $A$ is the path algebra of a finite quiver.
\end{enumerate}
\end{theorem}

Recall that an algebra $A$ is called \emph{Hopfian} if any $\kk$-algebra epimorphism from $A$ to itself is an automorphism. Moreover, we say $A$ is {\it strongly Hopfian} if $A[t_1,\cdots,t_n]$ is Hopfian for any $n\ge 0$. (see \cite[Definition 3.4]{LeWZ}). As a consequence, the following three large classes of strongly Hopfian algebras are Zariski cancellative once their centers are artinian. 

\begin{example}\cite[Lemma 3.5]{LeWZ}
The following algebras are strongly Hopfian.
\begin{enumerate}
    \item[(1)] Left or right noetherian algebras.
    \item[(2)] Locally finite $\mathbb N$-graded affine $\kk$-algebras.
    \item[(3)] Prime affine $\kk$-algebras satisfying a polynomial identity.  
\end{enumerate}
\end{example}

One may also consider the cancellation problem in the category of graded $\kk$-algebras. An isomorphism lemma is established for graded algebras in \cite{BZ2} stating that for any two connected graded $\kk$-algebras that are finitely generated in degree one, if they are isomorphic as ungraded algebras, then they must be isomorphic as graded algebras. As an application, we have the following result.

\begin{theorem}[\cite{BZ2}]
    Let $\Omega$ be the class of connected graded $\kk$-algebras finitely generated in degree one. If $A\in \Omega$ satisfies $Z(A)\cap A_1=\{0\},$ then $A$ is cancellative within $\Omega$, where $A_1$ is the set of elements of degree one in $A$.
\end{theorem}

The following results are proved in \cite{BR}, which indicate how the center of an algebra affects its cancellativity. The property of \emph{(strongly) retractable} and its weaker version of \emph{(strongly) detectable} are introduced in \cite{LeWZ} for non-domain noncommutative algebras parallel to the
discriminant computation method. Recall that \emph{the nilpotent radical} of a ring $R$ is defined to be the intersection of all prime ideals in $R$.

\begin{theorem}
   Let $R$ be a ring with center $Z$ and $N$ be the nilpotent radical of $R$. In each of the following cases, $R$ is strongly cancellative:
    \begin{enumerate}
        \item $Z/N$ is strongly retractable.
        \item $Z$ has Krull dimension zero,
        \item $Z$ is a finite direct sum of local rings.
    \end{enumerate}
\end{theorem}

For arbitrary affine $\kk$-domain of finite GK-dimension, \Cref{thm:LND} provides a general approach to the Zariski cancellation problem using ${\rm LND}$-rigidity. Discriminant computation so far is the most practical way to achieve that goal in the PI case.  

\begin{theorem}\cite[Theorem 5.2]{BZ1} 
Let $A$ be a PI $\kk$-domain. Suppose the $w$-discriminant of $A$ over its center is effective for some $w$. Then $A$ is strongly ${\rm LND}^{H}$-rigid. As a consequence, $A$ is strongly cancellative.
\end{theorem}

The above theorem can be applied to several important families of noncommutative PI-algebras. 

\begin{example}\cite[Theorem 0.8]{BZ1}
Let $A$ be a finite tensor product of the skew polynomial rings $\kk_{q}[x_{1}, \cdots, x_{n}]$ where $n$ is even and $q\in \kk\setminus\{0,1\}$ is a root of unity. Then $A$ is ${\rm LND}$-rigid. As a consequence, $A$ is cancellative.
\end{example}

In the example of skew polynomial rings $\kk_q[x_1,\cdots,x_n]$, we ask the following question when $n$ is odd. 

\begin{question}
\cite[Question 0.8]{TVZ}. Let $q\in \kk \setminus \{0, 1\}$ be a root of unity. Is the skew polynomial ring $\kk_{q}[x_{1},\cdots, x_{n}]$ cancellative for $n$ odd and $n\geq 3$?
\end{question}

Calculating the discriminants is usually difficult. Even if the discriminant is computable, it often lacks the dominating or effective property. Therefore, it is important to search for other means to prove the ${\rm LND}$-rigidity for a given algebra. Via the approach of {\it Nakayama automorphisms}, the cancellation problem is considered for certain special classes of non-PI AS-regular algebras of GK-dimension $3$ in \cite{LMZ}. In particular, many examples of cancellative algebras $A$ of GK-dimension no greater than 3 are given in \cite[Corollary 0.9] {LMZ}.

The cancellation problem is further considered for noetherian connected graded AS-regular algebras of GK-dimension $3$ or higher in \cite{TVZ}. The following more general results are proved in \cite{TVZ}.

\begin{theorem}\label{thm:0203}
\cite[Theorem 0.2 and 0.3]{TVZ} Suppose ${\rm char}(\kk)=0$. Let $A$ be a noetherian AS-regular algebra generated in degree one. If either of the following sets of hypotheses holds, then $A$ is cancellative:
\begin{enumerate}
\item ${\rm gl.dim}(A)=3$ and $A$ is not PI.
\item ${\rm GKdim}(Z(A))\leq 1$, and
 ${\rm gl.dim}(A/(t))=\infty$ for every homogeneous central element $t\in Z(A)$ of positive degree.
\end{enumerate}
\end{theorem}

The following conjecture is presented in \cite{TVZ}. 
\begin{conjecture}\label{conj:GKgldim}
\cite[Conjecture 0.4]{TVZ} Suppose ${\rm char}(\kk)=0$. Let $A$ be a noetherian finitely generated prime algebra.
\begin{enumerate}
\item [(1)] If ${\rm GKdim}(Z(A)) \leq 1$, then $A$ is cancellative.
\item [(2)] If ${\rm GKdim}(A)=3$ and $A$ is not PI, then $A$ is cancellative.
\end{enumerate}
\end{conjecture}

Regarding \Cref{thm:0203}, \Cref{conj:GKgldim} is just claiming that the global dimension hypotheses can be dropped if A is finitely generated prime or replaced by ${\rm GKdim}(A) = 3$ in part (2). This further leads to the following natural question from the perspective of the global or Krull dimension.

\begin{question}
\cite[Question 0.8]{LeWZ} Let $A$ be a $\kk$-algebra of global dimension one (respectively, Krull dimension
one). Is then $A$ cancellative?
\end{question}

It is important to mention that, under some natural assumptions, the PI quantized Weyl algebras are shown to be ${\rm LND}^H$-rigid and thus strongly cancellative using Poisson geometry \cite{LY}.

\section{Morita and Skew Cancellations}
It is natural to consider other generalizations of the Zariski cancellation problem. Indeed, one may consider the categorical versions of the ZCP. The notion of cancellation for abelian categories and derived categories are initiated in \cite{LuWZ}.  For any algebra $A$, let $M(A)$ denote the category of the right $A$-modules and $D(A)$ denote the derived category of right $A$-modules. The following definitions are proposed in \cite{LuWZ}.

\begin{definition}
Let $A$ be any algebra. 
\begin{itemize}
    \item[(1)]\cite[Definitions 0.1]{LuWZ} We say $A$ is {\it Morita cancellative} if the statement that $M(A[t])$ is equivalent to $M(B[t])$ for another algebra $B$ implies that $M(A)$ is equivalent to $M(B)$. 
    \item[(2)]\cite[Definitions 0.2]{LuWZ} We say $A$ is {\it derived cancellative} if the statement that 
$D(A[t])$ is triangulated equivalent to $D(B[t])$ for an algebra $B$
implies that $D(A)$ is triangulated equivalent to $D(B)$.
\end{itemize}
\end{definition}

From Morita theory, any two Morita equivalent algebras must have isomorphic centers. This suggests that many techniques used in the Zariski cancellation problem for non-domain noncommutative algebras based on their centers can be applied to the Morita cancellation problem. For example, according to \cite[Theorem 0.7]{LuWZ} a commutative domain is cancellative if and only if it is Morita cancellative if and only if $Z$ is derived cancellative. 

Moreover, a representation-theoretical version of discriminant was introduced in \cite{LuWZ} in terms of the module category which is well-behaved under Morita equivalence. As a consequence, relative ${\rm ML}_Z$-invariant and the notion of ${\rm LND}_Z$-rigidity with respect to the center $Z$ was extended to element-wise versions with respect to some representation-theoretical discriminant. Using these new ideas,  we have the Morita cancellation property for the following noncommutative algebras

\begin{theorem} 
The following algebras are Morita cancellative.
\begin{itemize}
\item[(1)] (\cite[Theorem 0.3]{LuWZ}) Every strongly Hopfian algebra with artinian center. 
    \item[(2)](\cite[Theorem 0.5]{LuWZ}) The path algebra $\kk Q$ for every finite quiver $Q$. 
    \item[(3)](\cite[Theorem 0.6]{LuWZ}) Every affine prime $\kk$-algebra of GK-dimension one over an algebraically closed field $\kk$. 
\item[(4)] (\cite[Corollary 0.8]{LuWZ}) Every non-PI Sklyanin algebra of global dimension three. 
\end{itemize}
\end{theorem}

The Morita cancellation is further studied in \cite{TZZ}, where the Morita version of the universally cancellative property is proposed in \cite[Definition 0.3]{TZZ} as a natural combination of being Morita cancellative and universally cancellative.  As a Morita version of \cite[Proposition 1.3]{BZ1}, one has the following result, which once again demonstrates the control of the center of $A$ over its cancellation property.
\begin{theorem} \cite[Theorem 0.4]{TZZ}
\label{xxthm0.4}
Let $A$ be an algebra with center $Z(A)$ being the base field
$\Bbbk$. Then $A$ is universally Morita cancellative.
\end{theorem}
As both a Morita version and a strengthened version of a partial combination of 
\cite[Theorem 4.1]{LeWZ} with \cite[Theorem 4.2]{LeWZ}, the following result is 
proved in \cite{TZZ}. Note that \emph{strongly Morita cancellation} can be defined analogously as strongly cancellation. We denote the nilradical of an algebra $A$ by $N(A)$.

\begin{theorem}\cite[Theorem 0.5]{TZZ}
\label{xxthm0.5}
Let $A$ be an algebra with center $Z$ such that either $Z$ or $Z/N(Z)$ is strongly
retractable {\rm{(}}respectively, strongly detectable{\rm{)}}, then $A$ is strongly Morita cancellative.
\end{theorem}
We remark that the hypotheses of being ``strongly Hopfian'' in the previous results such as \cite[Theorem 0.2, Theorem 4.2]{LeWZ}
and \cite[Theorem 0.3, Lemma 3.6, Theorem 4.2(2), Corollary 4.3, Corollary 7.3]{LuWZ} are not needed. So the above theorem provides a method to study the Morita Zariski cancellation problem for noncommutative algebras via their centers (e.g., \cite[Question 0.1]{LeWZ}). Recall that a commutative algebra is 
called {\it von Neumann regular} if it is reduced and has Krull dimension zero. 
One has the following result. 

\begin{corollary}\cite[Corollary 0.6]{TZZ}
\label{xxcor0.6}
Let $A$ be an algebra with center $Z$.
\begin{enumerate}
\item[(1)]
If $Z/N(Z)$ is generated by a set of units of $Z/N(Z)$,
then $Z$ and $A$ are strongly cancellative and strongly Morita
cancellative.
\item[(2)]
If $Z/N(Z)$ is a von Neumann regular algebra, then $Z$ and $A$ are strongly
cancellative and strongly Morita cancellative.
\item[(3)]
If $Z$ is a finite direct sum of local algebras, then $Z$ and $A$
are strongly cancellative and strongly Morita cancellative.
\end{enumerate}
\end{corollary}

The problem of skew cancellations was considered in \cite{AKP} and has been recently revisited in \cite{BHHV, Be}. Let $A$ be a $\kk$-algebra. Let $\sigma$ be a $\kk$-algebra automorphism of $A$ and $\delta$ be a $\sigma$-derivation of $A$. Then one can form the Ore extension, denoted by $A[t;\sigma,\delta]$, which shares many nice properties with the polynomial extension $A[t]$. An iterated Ore extension of an algebra $A$ is of the form
$$A[t_1;\sigma_1,\delta_1][t_2;\sigma_2, \delta]\cdots
[t_n;\sigma_n,\delta_n],$$
where $\sigma_{i}$ is an algebra automorphism of
$A_{i-1}:=A[t_1;\sigma_1,\delta_1]\cdots [t_{i-1};\sigma_{i-1},\delta_{i-1}]$
and $\delta_i$ is a $\sigma_i$-derivation of $A_{i-1}$. The reader is referred to \cite[Chapter 1]{MR} for more details, see also \cite{GoW}.

\begin{definition}\cite[Definition 0.7 and 0.8]{TZZ}
\label{xxdef0.7}
Let $A$ be an algebra.
\begin{enumerate}
\item[(1)]
We say $A$ is {\it skew cancellative} if any
isomorphism of algebras
\begin{equation}\label{eq:skewcan}
A[t;\sigma,\delta]\cong A'[t';\sigma',\delta']
\end{equation}
for another algebra $A'$, implies an isomorphism of algebras
$$A\cong A'.$$ 
\item[(2)]
We say $A$ is {\it $\sigma$-cancellative} if it is skew cancellative for Ore extensions \eqref{eq:skewcan} with $\delta=0$ and $\delta'=0$. 
\item[(3)]
We say $A$ is {\it $\delta$-cancellative} if it is skew cancellative for Ore extensions \eqref{eq:skewcan} with $\sigma={\rm Id}_{A}$ and $\sigma'={\rm Id}_{A'}$ .
\item[(4)]
We say $A$ is {\it $\sigma$-algebraically cancellative} if it is skew cancellative for Ore extensions \eqref{eq:skewcan} with locally algebraic $\sigma$ and $\sigma'$. 
\end{enumerate}
Moreover, we can define $A$ to be {\it strongly skew cancellative} if 
\[A[t_1;\sigma_1,\delta_1][t_2;\sigma_2, \delta_2]\cdots[t_n;\sigma_n,\delta_n]\cong A'[t'_1;\sigma'_1,\delta'_1][t'_2;\sigma'_2, \delta'_2]\cdots[t'_n;\sigma'_n,\delta'_n]\] 
implies that $A\cong A'$ for any $n\geq 1$.   
\end{definition}

Skew cancellation has appeared earlier in different notions. Given two Ore extensions $A[t; \delta_1]$ and $B[t; \delta_2]$, $A$ is said to be {\it Ore invariant} if the isomorphism $A[t;\delta_1] \cong B[t; \delta_2]$ implies $A\cong B$. When, moreover, the isomorphism between the Ore extensions carries $A$ onto $B$, then $A$ is said to be {\it strongly Ore invariant} \cite{AKP}. In \cite{AKP}, it is proved that certain abelian regular rings are strongly Ore invariant, and regular self-injective PI rings with no $\mathbb{Z}$-torsion are Ore invariant.

In \cite{Be}, it is proved that if the skew polynomial ring $R[t;\delta_1]$ and $\kk[x][t;\delta_2]$ are isomorphic and $\delta_2(x)\in \kk[x]$ has degree at least one, then $R\cong\kk[x]$. This result is generalized in \cite{BHHV} as follows.
\begin{theorem}\cite[Theorem 1.2]{BHHV}
Let $A$ be an affine commutative $\kk$-domain of Krull dimension one. Then we have the following.
\begin{itemize}
\item[(1)] $A$ is $\sigma$-cancellative.
\item[(2)] If $\kk$ has characteristic zero, then $A$ is $\delta$-cancellative.
\end{itemize}
\end{theorem}

However, the question of whether the cancellation property holds for any skew polynomial extensions of mixed type with coefficient rings being domains of Krull dimension one remains as an unresolved one in \cite{BHHV}.
\begin{question}
\cite[Question 5.7]{BHHV} Let $R$ be an affine commutative $\kk$-domain of Krull dimension one. Is any Ore extension $R[x; \sigma, \delta]$ skew cancellative?
\end{question}

Skew cancellation is further studied in \cite{TZZ} in terms of iterated Ore extensions and Ore extensions of special types. Moreover, an automorphism $\sigma$ of $A$ is called {\it locally algebraic} if every finite-dimensional subspace of $A$ is contained in a $\sigma$-stable finite dimensional subspace of $A$ \cite{Zh}. It is obvious that the identity map is locally algebraic. 

There are two new tools used in the skew cancellation problem. The first tool is the notion of divisor subalgebra, which was introduced in \cite{CYZ1} to solve Zariski cancellation and automorphism problems for certain noncommutative algebras. Let $A$ be a $\kk$-domain, and $F$ be any subset of $A$. We denote by $Sw(F)$ the set of $g\in A$ such that $f=agb$ for some $a, b\in A$ and
$0\neq f\in F$. That is, $Sw(F)$ is the set consisting of all the subwords of the
elements in $F$. Let us set $D_{0}(F)=F$ and inductively define $D_{n}(F)$ for $n\geq 1$ as the
$\Bbbk$-subalgebra of $A$ generated by $Sw(D_{n-1}(F))$. The
subalgebra ${\mathbb D}(F)=\bigcup_{n\geq 0} D_{n}(F)$ is called the {\it divisor
subalgebra} of $A$ generated by $F$. Any nonzero element $f\in A$ is called {\it controlling} if $\mathbb D(\{f\})=A$; see \cite[Examples]{TZZ} for controlling elements for various examples of noncommutative affine domains.

\begin{theorem}\cite[Theorem 0.9]{TZZ}
\label{xxthm0.9}
Let $A$ be an affine $\kk$-domain of finite GK-dimension. Suppose the unit element $1$ is controlling. Then $A$ is strongly $\sigma$-algebraically
cancellative. As a consequence, $A$ is strongly $\delta$-cancellative.
\end{theorem}

The second tool relies on a structural result of division algebras. Recall that a 
simple artinian ring $S$ is called {\it stratiform} over $\Bbbk$ \cite{Sc} if there is a chain 
of simple artinian rings $$S = S_n \supseteq S_{n-1} \supseteq \cdots \supseteq S_1 \supseteq S_0 =\Bbbk$$
where, for every $i$, either
\begin{enumerate}
\item[(i)]
$S_{i+1}$ is finite over $S_{i}$ on both sides; or
\item[(ii)]
$S_{i+1}$ is equal to the quotient ring of the
Ore extension $S_i[t_i; \sigma_i, \delta_i]$ for
an automorphism $\sigma_i$ of $S_i$ and $\sigma_i$-derivation
$\delta_i$ of $S_i$.
\end{enumerate}
Such a chain of simple artinian rings is called a stratification of
$S$. The {\it stratiform length} of $S$ is the number of steps in
the chain that are of type (ii). A crucial result proved in
\cite{Sc} is that the stratiform length is an invariant of $S$.
Moreover, a Goldie prime ring $A$ is called {\it stratiform} if its quotient
division ring, denoted by $Q(A)$, is stratiform. Many stratiform noetherian 
domains are proved to be skew cancellative. 

\begin{theorem}\cite[Theorem 0.10]{TZZ}
\label{xxthm0.10}
Let $A$ be a noetherian domain that is stratiform. Suppose
the unit element $1$ is controlling. Then $A$ is strongly skew cancellative
in the category of noetherian stratiform domains.
\end{theorem}

At the end of \cite{TZZ}, motivated by the results \cite[Lemma 4.3, Proposition 5.6]{BHHV}, it is established that every ${\rm LND}$-rigid algebra is $\delta$-cancellative. The following is the detailed statement. 
\begin{theorem}\cite[Theorem 5.4]{TZZ}
\label{xxthm5.4}
Suppose that $\Bbbk$ is a field of characteristic zero.
Let $A$ be an affine $\Bbbk$-domain of finite GK-dimension. Suppose
that ${\rm ML}(A)=A$. Then $A$ is $\delta$-cancellative.
\end{theorem}
We end this section by stating the following question from \cite{TZZ}.
\begin{question}
\cite[Question 5.6]{TZZ}. Suppose that $A$ is an affine domain of finite GK-dimension over a base field $\kk$ of characteristic zero. Suppose either ${\rm ML}(A)=A$ or $A$ is strongly ${\rm LND}$-rigid. Is $A$ strongly $\delta$-cancellative?
\end{question}

\section{Zariski cancellation for Poisson algebras}
Zariski cancellation problems can be asked for different types of associative and non-associative algebras. In this section, we focus on the Zariski cancellation problem for Poisson algebras. 

The notion of the Poisson bracket, first introduced by  Sim{\'e}on Denis Poisson, arises naturally in Hamiltonian mechanics and differential geometry. 

\begin{definition}
    A {\it Poisson algebra} is a commutative $\Bbbk$-algebra $A$ together with a bilinear form $\{-,-\}: A\times A \to A$ that is both a Lie bracket and a biderivation.  
\end{definition}
Poisson algebras have recently been studied intensively by many researchers with topics related to noncommutative discriminant \cite{BY,NTY}, representation theory of PI Sklyanin algebras  \cite{WWY1,WWY2}, Poisson Dixmier-Moeglin equivalences \cite{BLSM, LS, LWW}, Poisson enveloping algebras \cite{LOW, LWZ3, LWZ2, LWZ1}, invariant theory \cite{GVW} and so on. Moreover, the following question was asked in \cite{GaW} regarding the Zariski cancellation property.
\begin{question}[Zariski Cancellation Problem for Poisson Algebras]
When is a Poisson algebra $A$ cancellative? That is, does an isomorphism of Poisson algebras $A[t] \cong B[t]$ for another Poisson algebra $B$ imply an isomorphism $A \cong B$ as Poisson algebras?
\end{question}

We first restrict our attention to graded Poisson algebras. The main result below can be viewed as a Poisson version of \cite[Theorem 3.1]{BZ2}. Recall that any Poisson algebra $A$, the {\it Poisson center} of $A$ is denoted by $\mathcal Z_P(A):=\{a\in A\,|\, \{a,-\}=0\}$.   

\begin{theorem}\cite[Theorem 4.5]{GaW}
Let $A$ and $B$ be two connected graded Poisson algebras finitely generated in degree one. Suppose either $\mathcal Z_P(A)$ or $\mathcal Z_P(B)$ is generated in degree at least $2$. If $A[t_1,\dots,t_n]\cong B[t_1,\dots,t_n]$ as ungraded Poisson algebras for some $n\ge 1$, then $A\cong B$ as connected graded Poisson algebras.
\end{theorem}

The following result shows a deep connection between the Poisson center and Poisson cancellation.

\begin{theorem}
\label{thm.intro2}
Let $A$ be a Poisson algebra.
\begin{enumerate}
    \item [(a)] (\cite[Corollary 5.4]{GaW}) If $A$ is noetherian with artinian Poisson center, then $A$ is Poisson cancellative.
    \item [(b)] (\cite[Theorem 5.5]{GaW}) If $A$ has trivial Poisson center, then $A$ is Poisson cancellative.
\end{enumerate}
\end{theorem}

\noindent Theorem \ref{thm.intro2} (b) can be applied to show that Poisson integral domains of Krull dimension two which have nontrivial Poisson brackets are Poisson cancellative (\cite[Corollary 5.6]{GaW}). Here the non-triviality of the Poisson bracket plays an essential role since by \cite{Da, Fi} there are commutative domains of Krull dimension two that are not cancellative. 

It is important to point out that many theories developed in the (associative) noncommutative setting \cite{BZ1,BZ2} can be adapted to the Poisson setting. Hence the main work of \cite{GaW} is to modify the ideas of the noncommutative discriminant \cite{CPWZ1,CPWZ2} and the Makar-Limanov invariant \cite{Ma2} into the setting of Poisson algebras.

\begin{definition}
Let A be a Poisson algebra. A derivation $\delta$ of $A$ is called a {\it Poisson derivation} if 
\[
\delta(\{a,b\})~=~\{\delta(a),b\}+\{a,\delta(b)\} \text{ for all $a,b\in A$}.
\]
The {\it Poisson Makar-Limanov invariant} of $A$ is defined to be 
\[
{\rm PML}(A)~=~ \bigcap_{\delta\in {\rm PLND}(A)}\, {\rm ker}(\delta),
\]
 where ${\rm PLND}(A)$ denotes the space of all locally nilpotent Poisson derivations of A.
\end{definition}

\begin{theorem}\cite[Theorem 6.12]{GaW}
\label{thm.intro3}
Assume $\kk$ is a field of characteristic zero. Let $A$ be an affine Poisson domain over $\kk$ with finite Krull dimension. If $A$ has no nontrivial locally nilpotent Poisson derivations, then $A$ is Poisson cancellative.
\end{theorem}

In positive characteristic, we can prove a similar result by replacing Poisson derivations by higher Poisson derivations introduced in \cite{LL}.

It is important to mention that for Poisson algebras in characteristic zero, their Poisson centers are usually not large enough for us to emulate the definition of discriminant for noncommutative algebras by simply replacing the algebraic center with Poisson center. The following definition used the idea in \cite[\S 2]{LuWZ} to introduce the notion of Poisson discriminant from a representation-theoretic point of view.  

\begin{definition}\cite[Definition 7.1]{GaW}
\label{def.P}
Let $A$ be a Poisson algebra and let $\mathcal Z_P =\mathcal Z_P(A)$ be its Poisson center. Let $\mathcal P$ be a property defined for Poisson algebras that is invariant under isomorphism classes of Poisson algebras. We define the following terms for sets/ideals in $A$.
\begin{enumerate}
\item [(a)] ($\mathcal P$-locus) 
$L_{\mathcal P}(A) := \{ \fm \in {\rm Maxspec}(\mathcal Z_P) : A/\fm A \text{ has property $\mathcal P$}\}$.
\item [(b)] ($\mathcal P$-discriminant set) 
$D_{\mathcal P}(A) := {\rm Maxspec}(\mathcal Z_P) \backslash L_{\mathcal P}(A)$.
\item [(c)] ($\mathcal P$-discriminant ideal)
$I_{\mathcal P}(A) := 
%I(D_{\cP}) = 
\bigcap_{\mathfrak m \in D_{\mathcal P}(A)} \mathfrak m \subset \mathcal Z_P$.
\end{enumerate}
In the case that $I_{\mathcal P}(A)$ is a principal ideal, generated by $d \in \mathcal Z_P$, then $d$ is called the {\it $\mathcal P$-discriminant} of $A$, denoted by $d_{\mathcal P}(A)$. Observe that, if $\mathcal Z_P$ is a domain, $d_{\mathcal P}(A)$ is unique up to an element of $\mathcal Z_P^\times$.
\end{definition}

There is also a notion of effectiveness for Poisson discriminant. As in the noncommutative setting, effectiveness controls the locally nilpotent Poisson derivations and hence yields Poisson cancellation.  

\begin{theorem}\cite[Theorem 7.16]{GaW}
Let $A$ be an affine Poisson domain with affine Poisson center. If the Poisson discriminant exists and is effective either in $A$ or its Poisson center, then $A$ is Poisson cancellative.  
\end{theorem}

Next, we turn to the cancellation of quadratic Poisson polynomial algebras in three variables in comparison with their quantizations as AS-regular algebras of dimension three. By analyzing the Poisson center of an arbitrarily connected graded Poisson domain, it was shown that Poisson polynomial algebras in three variables with Poisson bracket either being quadratic or derived from a Lie algebra are cancellative.

\begin{theorem}
Let $\kk$ be a base field that is algebraically closed of characteristic zero. 
\begin{enumerate}
\item [(a)] (\cite[Theorem 3.11]{GaWY} and \cite[Corollary 3.14]{GaWY}) Let $A=\kk[x,y,z]$ be a quadratic Poisson algebra with nontrivial Poisson bracket. Then the $d$th Veronese Poisson subalgebra $A^{(d)}$ is Poisson cancellative for any $d\ge 1$. In particular, $A$ is Poisson cancellative. 
\item [(b)] (\cite[Theorem 3.15]{GaWY}) Let $\mathfrak g$ be a non-abelian Lie algebra of dimension $\le 3$. Then the symmetric algebra $S(\mathfrak g)$ on $\mathfrak g$ together with the Konstant-Kirillov bracket is Poisson cancellative. 
\end{enumerate}
\end{theorem}

As the Ore extensions of rings, there is a Poisson version of the extension of an arbitrary Poisson algebra, called Poisson-Ore extension, using a Poisson derivation $\alpha$ and a Poisson $\alpha$-derivation $\delta$ of the base Poisson algebra. In \cite{GaWY}, the \emph{skew Zariski cancellation problem}  was introduced for Poisson algebras in terms of any Poisson-Ore extension. That is, it asks whether or not an isomorphism between the Poisson-Ore extensions of two base Poisson algebras implies an isomorphism between the two base Poisson algebras. Like in the noncommutative algebraic setting \cite{BHHV,Be,TZZ}, similar invariants including the Poisson Makar-Limanov invariant, the divisor Poisson subalgebra, and the Poisson stratiform length can be defined for Poisson algebras, which play important roles in their skew cancellation problem. 

\begin{theorem}
Let $\kk$ be a base field of characteristic zero, and $A$ be a noetherian Poisson domain.
\begin{enumerate}
\item [(a)] (\cite[Theorem 4.4]{GaWY}) If $A$ is either Poisson simple of finite Krull dimension with units in  $\kk$ or affine of Krull dimension one, then $A$ is Poisson $\alpha$-cancellative.
\item [(b)] (\cite[Theorem 4.5]{GaWY})  If the Poisson Makar-Limanov invariant of $A$ equals $A$, then $A$ is Poisson $\delta$-cancellative. 
\item [(c)] (\cite[Theorem 4.8]{GaWY})  If $A$ has finite Krull dimension and the $1$-divisor Poisson subalgebra  of $A$ equals $A$, then $A$ is Poisson skew cancellative.
\item [(d)] (\cite[Theorem 4.15]{GaWY})
If $A$ is Poisson stratiform and the $1$-divisor Poisson subalgebra of $A$ equals $A$, then $A$ is Poisson skew cancellative in the category of noetherian Poisson stratiform domains.
\end{enumerate}
\end{theorem}

In the last part of this section, we will discuss some open questions regarding the Poisson cancellation problem. (see \cite[\S 8]{GaW} and \cite[\S 5]{GaWY} for more questions).

For any Poisson algebra A, there is a notion of Poisson universal enveloping algebra $U(A)$ (see \cite{Oh}), whose representation category is Morita equivalent to the category of Poisson modules over $A$. 

\begin{question}\cite[Question 8.7]{GaW}
Let $A$ be any affine Poisson algebra, and $U(A)$ be its Poisson universal enveloping algebra. What is the relationship between $U(A)$ being cancellative and A being Poisson cancellative?
\end{question}

In practice, many Poisson structures can be derived from the process of semiclassical limits of quantized coordinate rings, for instance, see \cite{Go}. 

\begin{question}\cite[Question 8.8]{GaW}
What is the relation between the cancellation property of a quantized coordinate ring and the Poisson cancellation property of its semiclassical limit?
\end{question}

Recall that, in \cite[Definition 0.2]{LuWZ}, a notion of {\it derived cancellation} is further introduced for any associative algebra by replacing the respective module category with its derived category. We are interested in its Poisson version. 

\begin{question}\cite[Question 5.4]{GaWY}
Can one define a suitable notion of \emph{derived Poisson cancellation}? Under what conditions does a Poisson algebra satisfy this condition?
\end{question}

\providecommand{\bysame}{\leavevmode\hbox to3em{\hrulefill}\thinspace}
\providecommand{\MR}{\relax\ifhmode\unskip\space\fi MR }
\providecommand{\MRhref}[2]{%

\href{http://www.ams.org/mathscinet-getitem?mr=#1}{#2} }
\providecommand{\href}[2]{#2}

\bibliography{Survey.bbl} 
\bibliographystyle{plain}
\end{document}